\documentclass{gtmon_a}
\pdfoutput=1

\proceedingstitle{Proceedings of the Nishida Fest (Kinosaki 2003)}
\conferencestart{28 July 2003}
\conferenceend{8 August 2003}
\conferencename{International Conference in Homotopy Theory}
\conferencelocation{Kinosaki, Japan}

\editor{Matthew Ando}
\givenname{Matthew}
\surname{Ando}

\editor{Norihiko Minami}
\givenname{Norihiko}
\surname{Minami}

\editor{Jack Morava}
\givenname{Jack}
\surname{Morava}

\editor{W Stephen Wilson}
\givenname{W Stephen}
\surname{Wilson}

\title[The Rothenberg--Steenrod spectral sequence]{On the
Rothenberg--Steenrod spectral sequence\\\vglue -3pt for the mod 3 cohomology of the
classifying space\\of the exceptional Lie group $E_8$}

\author{Masaki Kameko}
\givenname{Masaki}
\surname{Kameko}
\address{Toyama University of International Studies\\\newline
Toyama 930-1292\\
Japan}
\email{kameko@tuins.ac.jp}
\urladdr{}

\author{Mamoru Mimura}
\givenname{Mamoru}
\surname{Mimura}
\address{Department of Mathematics\\
Okayama University\\\newline
Okayama 700-8530\\
Japan}
\email{mimura@math.okayama-u.ac.jp}
\urladdr{}

\volumenumber{10}
\issuenumber{}
\publicationyear{2007}
\papernumber{12}
\startpage{213}
\endpage{226}

\doi{}
\MR{}
\Zbl{}

\arxivreference{}


\received{31 July 2004}
\revised{26 October 2005}
\accepted{}
\published{29 January 2007}
\publishedonline{29 January 2007}
\proposed{}
\seconded{}
\corresponding{}
\version{}

\makeatletter
\def\cnewtheorem#1[#2]#3{\newtheorem{#1}{#3}[section]
\expandafter\let\csname c@#1\endcsname\c@Theorem}

\AtBeginDocument{\let\bar\wbar\let\tilde\wtilde\let\hat\what}
\def\S{Section }

\newtheorem{Theorem}{Theorem}[section]
\cnewtheorem{Proposition}[Theorem]{Proposition}
\cnewtheorem{Lemma}[Theorem]{Lemma}

\makeatother  %  move after \newtheorem block

\makeautorefname{Theorem}{Theorem}
\makeautorefname{Proposition}{Proposition}
\makeautorefname{Lemma}{Lemma}

\newcommand{\MyMatrix}[3]{
\renewcommand{\arraystretch}{1.5}
\left(\begin{array}{c|c}
#1 & #2 \\
\hline 
0 & #3
\end{array}
\right)
}

\newcommand{\colsep}{\setlength{\arraycolsep}{.5mm}}

\begin{document}

\begin{htmlabstract} 
We show that the Rothenberg--Steenrod spectral sequence converging to the
mod 3 cohomology of the classifying space of the exceptional Lie group
E<sub>8</sub> does not collapse at the E<sub>2</sub>&ndash;level.
\end{htmlabstract}

\begin{abstract} 
We show that the Rothenberg--Steenrod spectral sequence converging to the 
mod 3 cohomology of the classifying space of the exceptional Lie group $E_8$ 
does not collapse at the $E_2$--level.
\end{abstract}

\maketitle

\section{Introduction}
\label{intro}
One of the most powerful tools in the study of the mod $p$ 
cohomology of classifying spaces of connected compact Lie groups 
is the Rothenberg--Steenrod spectral sequence.
For a connected compact Lie group $G$, there is a 
strongly convergent first quadrant spectral sequence of graded $\mathbb{F}_p$--algebras
\[\{ E^{p,q}_r, d_r\co E_r^{p,q}\to E_{r}^{p+r, q-r+1}\}\]
such that 
$\displaystyle
E_2^{p,q} = \mathrm{Cotor }_{H^{*}\!G}^{\,p,q}(\mathbb{F}_p, \mathbb{F}_p)
$
and $
\displaystyle E_{\infty}=\mathrm{gr}\ H^{*}\!BG$.
The Rothenberg--Steenrod spectral sequence, sometimes mentioned as the 
 Eilenberg--Moore spectral sequence, 
has been successful in computing the mod $p$ cohomology
of classifying spaces of  connected compact Lie groups.
In all known cases, for an odd prime $p$, the Rothenberg--Steenrod spectral sequence
converging to the mod $p$ cohomology of the classifying space of a connected 
compact Lie group collapses at the $E_2$--level.
So, one might expect that
this collapse should always occur.
%for an odd prime $p$,
%the Rothenberg--Steenrod spectral sequence converging to the mod $p$ cohomology
%of the classifying space of a connected compact Lie group always 
%collapses at the $E_2$--level.
In this paper, however, we 
show that this is not the case for the mod $3$ cohomology of the classifying space of the exceptional Lie group $E_8$.

\begin{Theorem} \label{noncollapsing}
The Rothenberg--Steenrod spectral sequence converging to 
the mod $3$ cohomology of the classifying space of the exceptional Lie group $E_8$ does not collapse at the $E_2$--level.
\end{Theorem}

We prove \fullref{noncollapsing} by computing the ring of invariants of
the mod $3$ cohomology of a nontoral elementary abelian $3$--subgroup of $E_8$. 
According to Andersen, Grodal, M\o ller and Viruel \cite{Andersen},
up to conjugates, there are exactly two maximal nontoral elementary abelian $3$--subgroups, which they call $\smash{E_{E_8}^{5a}}$ and $\smash{E_{E_8}^{5b}}$. 
They described the action of Weyl groups on these 
nontoral elementary abelian 
$3$--subgroups explicitly. 
In this paper, we compute the ring of invariants 
of the polynomial part of the mod $3$ cohomology of $BE_{E_8}^{5a}$.
By comparing degrees of algebra generators of the above ring of invariants with those of algebra generators of 
the cotorsion product
$\mathrm{Cotor}_{H^{*}\!E_8}(\mathbb{F}_3, \mathbb{F}_3)$ computed by Mimura and Sambe \cite{MimuraSambe},
we prove \fullref{noncollapsing}. 

In \fullref{Sec2}, we set up 
a tool, \fullref{main}, for the computation of certain rings of invariants.
In \fullref{Sec3}, we recall some facts on maximal nontoral elementary abelian $p$--subgroups of simply connected compact simple Lie groups
and their Weyl groups.  Then using \fullref{main}, we compute 
some of rings of invariants of the above Weyl groups.
In \fullref{Sec4},  we complete the proof of \fullref{noncollapsing}.

The first named author was partially supported by the Japan Society for the Promotion of Science, Grant-in-Aid for Scientific Research (C) 13640090.

\section{Invariant theory}
\label{Sec2}

In this section, we consider the invariant theory 
over the finite field $\mathbb{F}_q$ of $q$ elements where $q=p^{k}$ with $k\geq 1$  and $p$ is a prime number.
For a finite set $\{v_1, \ldots, v_n\}$, we denote by $\mathbb{F}_q\{v_1, \ldots, v_n\}$ the $n$--dimensional
vector space over $\mathbb{F}_q$ spanned by $\{v_1, \ldots, v_n\}$.

Let us write $GL_{n}(\mathbb{F}_q)$ for
the set of invertible $n\times n$ matrices whose entries are in 
$\mathbb{F}_q$.
We also write $M_{m,n}(\mathbb{F}_q)$ for the set of $m\times n$ matrices whose entries  are in $\mathbb{F}_q$.
Let $G$ be a subgroup of $GL_{n}(\mathbb{F}_q)$.
The group $G$ acts on 
the $n$--dimensional vector space
$V=\mathbb{F}_q\{ v_1, \ldots, v_n\}$ as follows:
for $g$ in  $G$, 
\[
g v_i=\sum_{j=1}^{n} a_{j,i}(g)v_{j},
\]
where $a_{i,j}(g)$ is the $(i,j)$--entry of the matrix $g$.
We denote by 
$\{x_1, \ldots, x_{n}\}$ 
the dual basis of $\{v_1, \ldots, v_n\}$ and write $V^{*}$
 for the dual of $V$, that is,
 \[
 V^{*}=\mathrm{Hom}_{\mathbb{F}_q}(V, \mathbb{F}_q)=\mathbb{F}_q\{x_1, \ldots, x_n\}.
 \]
$$\mathbb{F}_q[V]=\mathbb{F}_q[x_1, \ldots, x_n]\leqno{\hbox{We denote by}}$$ the polynomial algebra over $\mathbb{F}_q$ in $n$ variables $x_1$, \ldots, $x_n$.
Then the group $G$ acts on both $V^*$ 
and $\mathbb{F}_{q}[V]$ as follows:
for $g$ in $G$, 
\begin{align*}
(gx)(v)&=x (g^{-1}v) &&\text{for } x \text{ in } V^*, v \text{ in } V;&
\\
g(y\cdot z)&=(gy)\cdot (gz) &&\text{for } y, z \text{ in } \mathbb{F}_q[V].&
\end{align*}
Using entries of a matrix $g\in G$, we may describe the action of $g$ as follows:
\begin{Proposition}\label{action1}
\[
gx_{i}=\sum_{j=1}^{n}a_{i,j}(g^{-1})x_j.
\]
\end{Proposition}
\begin{proof}
\[
(g x_{i})(v_{j})=x_{i}(g^{-1}v_{j})
=x_{i}\Big(\sum_{k=1}^{n} a_{k,j}(g^{-1})v_{k}\Big)
=a_{i, j}(g^{-1}).
\proved
\]
\end{proof}

In order to prove \fullref{main} below, 
we recall a strategy of Wilkerson \hbox{\cite[\S3]{Wilkerson}}.
It can be stated in the following form.
\begin{Theorem} \label{wilkerson}
Let $G$ be a subgroup of $GL_{n}(\mathbb{F}_q)$ acting on $V$
as above.
Let  $f_1, \ldots, f_n$ be homogeneous polynomials in 
$\mathbb{F}_q[V]$.
We have
\[\mathbb{F}_q[V]^{G}=\mathbb{F}_q[f_1, \ldots, f_n] \]
if and only if  the following three conditions hold:
\begin{enumerate}
\item  $f_1, \ldots, f_n$ are $G$--invariant;
\item  
$\mathbb{F}_q[V]$ is integral over the
subalgebra $R$ generated by 
$f_1, \ldots, f_n$;
\item  $\deg f_1 \cdots \deg f_n=|G|$.
\end{enumerate}
\end{Theorem}
In the statement of \fullref{wilkerson},
$\deg f$ is the homogeneous degree of $f$, that is, we define
the degree $\deg x_i$ of indeterminate $x_i$ to be $1$.
For the proof of this theorem, 
we refer the reader to Corollaries 2.3.2 and 5.5.4 and Proposition 5.5.5
in Smith's book \cite{Smith} and Wilkerson's paper \cite[\S3]{Wilkerson}.

To state \fullref{main},  we need to set up the notation.
Let $G_1 \subset GL_{m}(\mathbb{F}_q)$ and 
$G_2 \subset GL_{n-m}(\mathbb{F}_q)$. 
Let $V_1=\mathbb{F}_q\{v_1, \ldots, v_m\}$ and $V_2=\mathbb{F}_q\{v_{m+1}, \ldots, v_n\}$.
Let $G_1$ and $G_2$ act on $V_1$ and $V_2$ by
\[
g_1 v_i=\sum_{k=1}^{m} a_{k,i}(g_1) v_{k} \quad\mbox{and}\quad g_2 v_{j}=\sum_{k=1}^{n-m} a_{k, j-m}(g_2) v_{m+k},
\]
respectively, where $i=1, \ldots, m$ and $j=m+1, \ldots, n$.
The following proposition is immediate from the definition and \fullref{action1}.

\begin{Proposition}\label{action2} The following hold:
\begin{enumerate}
\item If $f(x_1, \ldots, x_m) \in \mathbb{F}_q[V_1]$ is $G_1$--invariant, then for all $g_1\in G_1$, we have
\[
f\Big(\sum_{k=1}^{m} a_{1, k}(g_1^{-1})x_{k}, \ldots, \sum_{k=1}^{m}a_{m,k}(g_1^{-1})x_k\Big)=f(x_1, \ldots, x_m);
\]
\item If
$f(x_{m+1}, \ldots, x_{n})\in \mathbb{F}_q[V_2]$ is $G_2$--invariant, then for all $g_2\in G_2$, we have
\[
f\Big(\sum_{k=1}^{n-m} a_{1, k}(g_2^{-1})x_{m+k}, \ldots, \sum_{k=1}^{n-m}a_{n-m,k}(g_2^{-1})x_{m+k}\Big)
=f(x_{m+1}, \ldots, x_n).
\]
\end{enumerate}
\end{Proposition}

Suppose that $G$ consists of  the matrices of the form
\[
\MyMatrix{g_1}{m_0}{g_2},
\]
where $g_1\in G_1\subset GL_{m}(\mathbb{F}_q)$,
$g_2\in G_2 \subset GL_{n-m}(\mathbb{F}_q)$
and $m_0\in M_{m, n-m}(\mathbb{F}_q)$.
We denote respectively by
 $\bar{G}_0$, $\bar{G}_1$, $\bar{G}_2$  the subgroups of $G$ 
consisting of matrices of the form
\[
\MyMatrix{1_{m}}{m_0}{1_{n-m}}, \MyMatrix{g_1}{0}{1_{n-m}}, \MyMatrix{1_{m}}{0}{g_2},
\]
where $g_1\in G_1$, $g_2\in G_2$, $m_0\in M_{m,n-m}(\mathbb{F}_q)$ and $1_k$ is the identity matrix in $M_{k,k}(\mathbb{F}_q)$.
We denote by $\bar{g}_1$, $\bar{g}_2$ the elements in $\bar{G}_1$, $\bar{G}_2$  corresponding to $g_1$, $g_2$, respectively.

Considering $V_2^{*}$ as a subspace of $V^{*}$,  let us define $\mathcal{O}X$ in $\mathbb{F}_q[V][X]$ by
\[
\mathcal{O}X=\prod_{x\in V_2^{*}}(X+x).
\]
The following proposition is well-known (see Wilkerson \cite[\S1]{Wilkerson}). 

\begin{Proposition}\label{odickson}
There are $c_{n-m, k}\in \mathbb{F}_q[V_2]$ such that
\[
\mathcal{O}X = \sum_{k=0}^{n-m} (-1)^{n-m-k}c_{n-m,k}X^{q^k},
\]
where $c_{n-m, n-m}=1$.
\end{Proposition}

Now, we state our main theorem of this section.

\begin{Theorem}\label{main}
With the above assumption on $G$, 
suppose that rings of invariants 
$
\mathbb{F}_q[V_1]^{G_1}
$
and 
$\mathbb{F}_q[V_2]^{G_2}
$
are polynomial algebras
\[ 
\mathbb{F}_q[f_1, \ldots, f_m] \quad\mbox{and} \quad
\mathbb{F}_q[f_{m+1}, \ldots, f_{n}],
\]
respectively, where $f_1, \ldots, f_m$ are homogeneous polynomials in $m$ variables $x_1, \ldots, x_m$ and
$f_{m+1}, \ldots, f_n$ are homogeneous polynomials in $(n-m)$ variables $x_{m+1}, \ldots, x_{n}$.
Then the ring of invariants
$
\mathbb{F}_q[V]^{G}$ 
is also a polynomial algebra 
\[
\mathbb{F}_q[\bar{f}_1, \ldots, \bar{f}_{m}, f_{m+1}, \ldots, f_{n}],
\]
where  for $i=1$, \ldots, $m$,
\[
\bar{f}_i=f_{i}(\mathcal{O}x_1, \ldots, \mathcal{O}x_m).
\]
\end{Theorem}

To prove \fullref{main}, 
we verify that the conditions (1), (2) and (3) in \fullref{wilkerson} hold for
the polynomials $\bar{f}_1, \ldots, \bar{f}_m$, 
$f_{m+1}, \ldots, f_n$ in \fullref{main}.

{\bf Step 1}\qua To prove that $\bar{f}_1, \ldots, \bar{f}_m, f_{m+1}, \ldots, f_n$ are $G$--invariant, 
it suffices to prove the following propositions.

\begin{Proposition} \label{invariant1}
Suppose $f(x_1, \ldots, x_m)\in \mathbb{F}_q[V_1]^{G_1}$.
Then $f(\mathcal{O}x_1, \ldots, \mathcal{O}x_m)$ is $G$--invariant in $\mathbb{F}_q[V]$.
\end{Proposition}

\begin{Proposition}\label{invariant2}
Suppose $f(x_{m+1}, \ldots, x_{n})\in \mathbb{F}_q[V_2]^{G_2}$. 
Then $f(x_{m+1}, \ldots, x_{n})$ is $G$--invariant in $\mathbb{F}_q[V]$.
\end{Proposition}

To simplify the argument, we use the following.

\begin{Lemma} \label{o0}
$f\in \mathbb{F}_q[V]$ is $G$--invariant
if $f$ is $\bar{G}_0$--invariant, $\bar{G}_1$--invariant and $\bar{G}_2$--invariant.
\end{Lemma}
\begin{proof}
We may express each $g$ in $G$ as a product of 
elements $\bar{g}_0$, $\bar{g}_1$, $\bar{g}_2$ 
in $\bar{G}_0$, $\bar{G}_1$, $\bar{G}_2$ respectively, say $g=\bar{g}_0\bar{g}_1\bar{g}_2$ as follows:
\[
\MyMatrix{g_1}{m_0}{g_2}=\MyMatrix{1_m}{m_0g_{2}^{-1}}{1_{n-m}}
\MyMatrix{g_1}{0}{1_{n-m}}\MyMatrix{1_m}{0}{g_2}.\proved
\]
\end{proof}

Firstly, we prove \fullref{invariant1}.

\begin{Lemma}\label{olinear}
$\mathcal{O}$ is an $\mathbb{F}_q$--linear homomorphism from 
$V^{*}$ to 
$\mathbb{F}_q[V]$.
\end{Lemma}
\begin{proof}
For $\alpha$, $\beta\in \mathbb{F}_q$ and for $x$, $y\in V^*$, 
we have $(\alpha x+ \beta y)^{q^{k}}=\alpha x^{q^k}+\beta y^{q^k}\!$.
By \fullref{odickson}, we have
\[
\eqalignbot{\mathcal{O}(\alpha x+\beta y)&= \sum_{k=0}^{n-m}
(-1)^{n-m-k}c_{n-m,k} (\alpha x+\beta y)^{q^k}\cr
&= \sum_{k=0}^{n-m} (-1)^{n-m-k} c_{n-m, k}(\alpha x^{q^k}+\beta y^{q^k})\cr
&=\alpha \mathcal{O}x+\beta \mathcal{O}y.
}
\proved
\]
\end{proof}

\begin{Lemma} \label{o1}
The following hold for $k=1, \ldots, m$:
\begin{enumerate}
\item $\displaystyle \bar{g}_0 \mathcal{O}x_k=\mathcal{O}x_{k}$; 
\item $\displaystyle \bar{g}_1\mathcal{O}x_{k}=\sum_{\ell=1}^{m} a_{k,\ell}(g_{1}^{-1}) \mathcal{O}x_{\ell}$;
\item $\displaystyle \bar{g}_2\mathcal{O}x_{k}=\mathcal{O}x_k$.
\end{enumerate}
\end{Lemma}
\begin{proof}
(1)\qua We have $\bar{g}_0x_k=x_k+y$ for some $y$ in $\smash{V_2^{*}}$, and $\bar{g}_0x=x$ for any $x$ in $V_2^*$.  Then
$y+x$ ranges over $\smash{V_2^*}$ as $x$ ranges over $\smash{V_2^*}$. Hence, we have
\[
\bar{g}_0 \mathcal{O}x_k 
= \displaystyle \prod_{x \in V_2^*} (x_k+y+x) 
= \mathcal{O}x_k.
\]
(2)\qua We have $\bar{g}_1x_{k}=\sum_{\ell=1}^{m} a_{k,\ell}(g_{1}^{-1})x_{\ell}$, and 
$\bar{g}_1x=x$ for any $x$ in $V_2^*$. 
Hence, by \fullref{olinear}, we have
\[
\colsep
\begin{array}{rcl}
\bar{g}_1 \mathcal{O}x_k & = &\displaystyle \prod_{x \in V_2^*} \Big(
\sum_{\ell=1}^{m} a_{k,\ell}(g_1^{-1})x_\ell+x\Big) \\
&=& \displaystyle \mathcal{O}\Big(\sum_{\ell=1}^{m} a_{k,\ell}(g_1^{-1})x_{\ell}\Big) \\
&=& \displaystyle \sum_{\ell=1}^{m} a_{k, \ell}(g_{1}^{-1})\mathcal{O}x_{\ell}.
\end{array}
\]
(3)\qua We have $\bar{g}_2x_{k}=x_{k}$. Then $\bar{g}_2x$ ranges over $V_2^*$ as $x$ ranges over $V_2^*$.
Hence,
\[
\bar{g}_2 \mathcal{O}x_{k} = \displaystyle \prod_{x \in V_2^*}(x_{k}+\bar{g}_2x) = \mathcal{O}x_k.
\proved
\]
\end{proof}
\begin{proof}[Proof of \fullref{invariant1}]
By \fullref{o0}, it suffices to show that $$f(\mathcal{O}x_1, \ldots, \mathcal{O}x_m)$$ is $\bar{G}_i$--invariant 
 for $i=0$, $1$, $2$.
By \fullref{o1} (1) and (3), it is clear that the above element is invariant with respect to the action of
$\bar{G}_i$ for $i=0$, $2$.
By \fullref{o1} (2) and by \fullref{action2} (1), we have
\[
\colsep
\begin{array}[b]{rcl}
\bar{g}_1f(\mathcal{O}x_1, \ldots, \mathcal{O}x_m)&=&\displaystyle
f\Big(\sum_{k=1}^{m}a_{1,k}(g_1^{-1})\mathcal{O}x_{k},\ldots,
\sum_{k=1}^{m} a_{m,k}(g_1^{-1})\mathcal{O}x_{k}\Big)\\
&=&f(\mathcal{O}x_1, \ldots, \mathcal{O}x_m).
\end{array}
\proved
\]
\end{proof}

Secondly, we prove \fullref{invariant2}.

\begin{Lemma} \label{o2} The following hold for $k=m+1, \ldots, n$:
\begin{enumerate}
\item $\bar{g}_0x_{k}=x_{k}$;
\item $\bar{g}_1x_{k}=x_{k}$;
\item $\displaystyle \bar{g}_2x_{k}=\sum_{\ell=1}^{n-m} a_{k-m, \ell}(g_2^{-1}) x_{m+\ell}$.
\end{enumerate}
\end{Lemma}
\begin{proof}
(1) and (2) are immediate from the definitions of $\bar{g}_0$ and $\bar{g}_1$. 
 (3) follows immediately from the fact that 
\[
a_{k, \ell}(\bar{g}_{2}^{-1})=a_{k-m, \ell-m}(g_{2}^{-1})\]
 for $\ell\geq m+1$ and that
$a_{k, \ell}(\bar{g}_2^{-1})=0$ for $\ell \leq m$.
\end{proof}
\begin{proof}[Proof of \fullref{invariant2}]
As in the proof of \fullref{invariant1}, 
it suffices to show that
$$f(x_{m+1}, \ldots, x_n)$$
is $\bar{G}_i$--invariant
for $i=0$, $1$, $2$.
It is clear from \fullref{o2} (1)--(2)  that the above element is $\bar{G}_i$--invariant
for $i=0$, $1$.
By \fullref{o2} (3) and by \fullref{action2} (2), 
\[
\colsep
\begin{array}[b]{rcl}
\bar{g}_{2}f(x_{m+1}, \ldots, x_{n})
&=&\displaystyle
f\Big(\sum_{k=1}^{n-m} a_{1,k}(g_2^{-1})x_{m+k}, \ldots, 
\sum_{k=1}^{n-m} a_{n-m,k}(g_2^{-1})x_{m+k}\Big)\\
&=& f(x_{m+1}, \ldots, x_{n}).
\end{array}
\proved
\]
\end{proof}

{\bf Step 2}\qua We prove that the inclusion $R\to \smash{\mathbb{F}_q[V]}$ is an integral extension,
for $R$ the subalgebra of $\mathbb{F}_q[V]$ 
generated by $\smash{\bar{f}_1, \ldots, \bar{f}_m}$, $f_{m+1}, \ldots, f_{n}$.
Let $S$ be the subalgebra of $\mathbb{F}_q[V]$
generated by $\smash{\bar{f}_1, \ldots, \bar{f}_m}$, $c_{n-m, 0}, \ldots, c_{n-m,n-m-1}$.
Since $G_2\subset GL_{n-m}(\mathbb{F}_q)$, we see that $c_{n-m, k}\in R$. So, 
$S$ is a subalgebra of $R$.
Therefore, it suffices to prove the following proposition.

\begin{Proposition}
For $k=1, \ldots, n$, the element $x_k$ is integral over $S$.
\end{Proposition}

\begin{proof}
Firstly, we prove that $x_k$ is integral over $S$ for $k=1, \ldots, m$.
By \hbox{\fullref{wilkerson}}, $x_k$ is integral over $\mathbb{F}_q[V_1]^{G_1}$.
Hence, there exists a monic polynomial
$F(X)$ and polynomials $\varphi_{j}$'s over $\mathbb{F}_q$ in $m$ variables 
for $j=0, \ldots, r-1$ such that
\[
F(X) = X^r +\sum_{j=0}^{r-1} \varphi_{j}(f_1(x_1, \ldots, x_m), \ldots, f_m(x_1, \ldots, x_m)) X^{j}
\]
and that
 $F(x_k)=0$ in $\mathbb{F}_q[x_1, \ldots, x_m]$.

Replacing $x_i$ in the equality $F(x_k)=0$ above by $\mathcal{O}x_i$ for $i=1, \ldots, m$, 
 we have the following equality in $\mathbb{F}_q[V]$:
 \begin{gather*}
 (\mathcal{O}x_k)^r +\sum_{j=0}^{r-1} \varphi_{j}(f_1(\mathcal{O}x_1, \ldots, \mathcal{O}x_m), \ldots, 
f_m(\mathcal{O}x_1, \ldots, \mathcal{O}x_m)) (\mathcal{O}x_k)^{j}=0.
\\
F'(X) =(\mathcal{O}X)^{r}+\sum_{j=0}^{r-1}
\varphi_{j}(\bar{f}_1, \ldots, \bar{f}_m)(\mathcal{O}X)^{j}.
 \tag*{\hbox{Let}}\end{gather*}
By \fullref{odickson},  $F'(X)$ is a monic polynomial in $S[X]$. 
Since, by definition, $\bar{f}_i=f_i(\mathcal{O}x_1, \ldots, \mathcal{O}x_m)$,   
it is clear that $F'(x_k)=0$ in $\mathbb{F}_q[V]$. 
 Hence $x_k$ is integral over $S$. 
 
Secondly, we verify that $x_k$ is integral over $S$ for $k=m+1, \ldots, n$.
 By \fullref{odickson}, $\mathcal{O}X$ is a monic polynomial in $S[X]$.
 It is immediate from the definition that $\mathcal{O}x=0$ for $x\in V_2^{*}$.
 Therefore,  $x_k$ is integral over $S$. 
  \end{proof}
 
{\bf Step 3}\qua Finally, we compute the product of degrees of $\bar{f}_1, \ldots, \bar{f}_m$, $f_{m+1}, \ldots, f_{n}$.
 Since $\deg \mathcal{O}x$ is of degree $q^{n-m}$  for $x \in V^{*}$, we have
 $$\deg \bar{f}_i=\deg f_i \cdot q^{n-m}. $$ 
 By \fullref{wilkerson},  we have 
 \begin{align*}\deg f_1 \cdots \deg f_m=|&G_1|\\
 \deg f_{m+1}\cdots \deg f_{n}=&|G_2|. \tag*{\hbox{and}}\\
 %\end{gather*}
% Therefore,
% \[
%\colsep
% \begin{array}{rcl}
 \tag*{\hbox{Therefore}} \deg \bar{f}_1 \ \cdots \deg \bar{f}_m\cdot  \deg f_{m+1} \cdots  \deg f_{n}
 &= \deg f_1 \cdots  \deg f_{n}\cdot q^{m(n-m)}\\
 &= |G_1| \cdot |G_2| \cdot q^{m(n-m)} \\
 &= |G|.
%& \end{array}
% \]
\end{align*}
 This completes the proof of \fullref{main}.

\section{Rings of invariants of Weyl groups}
\label{Sec3}
Let $p$ be an odd prime. Let $G$ be a compact Lie group. 
We write $H_*BG$ and $H^*BG$ for the mod $p$ homology and cohomology of the classifying space $BG$ of $G$.
We write $A$ for an elementary abelian $p$--subgroup of the compact Lie group $G$.
Let 
\[
\Gamma H^{*}\!BG=H^{*}\!BG/\sqrt{0}, 
\]
where $\sqrt{0}$ is the ideal of nilpotent elements in $H^*BG$.
It is clear that $\Gamma H^{*}\!BA$ is a polynomial algebra 
\[
\Gamma H^*\!BA=\mathbb{F}_p[t_1, \ldots, t_n],
\]
where the cohomological degree of each $t_i$ is 2 and $n$ is the rank of $A$. 
We called it the polynomial part of $H^*BA$ in \fullref{intro}.

Choosing a basis for $A$, we may consider the action of  $GL_{n}(\mathbb{F}_p)$  on $A$.
We recall the relation between the action of  $GL_{n}(\mathbb{F}_p)$ on $A$
and the one on $\Gamma H^*BA$.
For the sake of notational simplicity, let $V=H_1BA$.
On the one hand,  $V$ is
identified with $A$ as a $GL_{n}(\mathbb{F}_p)$--module, where $g\in GL_{n}(\mathbb{F}_p)$ acts on 
$V$ as the induced homomorphism $Bg_{*}$.
As an $\mathbb{F}_p$--algebra, $\Gamma H^{*}\!BA$ is isomorphic to
$\mathbb{F}_p[V]$.
As in the previous section, we may 
consider the $GL_{n}(\mathbb{F}_p)$--module structure 
on $\mathbb{F}_p[V]$.
On the other hand,  
 $\smash{GL_{n}(\mathbb{F}_p)}$ acts on 
$\smash{H^{*}\!BA}$ by 
$g x = \smash{B(g^{-1})^{*} x}$, where $\smash{x\in H^*BA}$ and $g\in GL_{n}(\mathbb{F}_p)$. 
The relation between these actions is given by the following proposition.

\begin{Proposition} As a $GL_{n}(\mathbb{F}_p)$--module,
 $\Gamma H^{*}\!BA=\mathbb{F}_p[V]$.
\end{Proposition}

\begin{proof}
The Bockstein homomorphism induces an isomorphism of 
$GL_{n}(\mathbb{F}_p)$--modules 
\[
\beta\co H^1BA\to \Gamma H^{2}BA.
\]
Since, for $x \in V^*=H^1BA$, $v\in V=H_1BA$, we have
\[
(g x )(v) = (x) (g^{-1}v) 
= (x) ({B(g^{-1})}_*v) 
= ({B(g^{-1})}^{*}x)(v),
\]
we see that 
$\Gamma H^2BA=H^1BA=V^*$ as $GL_{n}(\mathbb{F}_p)$--modules.
Hence, we may conclude that $\Gamma H^*BA=\mathbb{F}_p[V]$ as 
$GL_{n}(\mathbb{F}_p)$--modules.
\end{proof}

The Weyl group $W(A)=N_{G}(A)/C_{G}(A)$ acts on 
$A$ as inner automorphisms.
So, we have the action of $W(A)$ on $\Gamma H^*BA$. 
Choosing a basis for $A$, we consider the Weyl group 
$W(A)$ as a subgroup of $GL_{n}(\mathbb{F}_p)$. 

In this section, we compute rings of invariants of Weyl groups of the polynomial part of
the mod $p$ cohomology of the classifying spaces of maximal nontoral elementary abelian $p$--subgroups
of simply connected compact simple Lie groups.

It is well-known that for an odd prime $p$, a simply connected compact simple Lie group $G$ does not have nontoral elementary abelian $p$--subgroups 
except for the cases $p=5$,  $G=E_8$, and $p=3$, $G=F_4$, $E_6$, $E_7$, $E_8$.
Andersen, Grodal, M\o ller and Viruel \cite{Andersen} described Weyl groups of maximal 
 nontoral elementary abelian $p$--subgroups
and their action on the underlying elementary abelian $p$--subgroup explicitly for $p=3$, $G=E_6$, $E_7$, $E_8$.
Up to conjugate, there are only 6 maximal nontoral elementary abelian $p$--subgroups of simply connected compact simple Lie groups.
For $p=5$, $G=E_8$ and for $p=3$, $G=F_4$, $E_6$, $E_7$, 
there is one maximal nontoral elementary abelian $p$--subgroup for each $G$. We call them
$\smash{E_{E_8}^3}$, $\smash{E_{F_4}^3}$, $\smash{E_{3E_6}^4}$, $\smash{E_{2E_7}^4}$, following the notation in \cite{Andersen}.
For $p=3$, $G=E_8$, there are two maximal nontoral elementary abelian $p$--subgroups, say $\smash{E_{E_8}^{5a}}$ and $\smash{E_{E_8}^{5b}}$,
where the superscript indicates the rank of elementary abelian $p$--subgroup.
For a detailed  account on nontoral elementary abelian $p$--subgroups, 
we refer the reader to \cite[\S8]{Andersen} and its references.

In this section, we compute
\[
(\Gamma H^{*}\!BA)^{W(A)}
\]
for $A=E_{3E_6}^4$, $E_{2E_7}^4$, $E_{E_8}^{5a}$ using \fullref{main}.  

\begin{Proposition}\label{ring} 
We have the following isomorphisms of graded $\mathbb{F}_p$--algebras:
\begin{enumerate}
\item For $p=5$, $G=E_8$, $A=E_{E_8}^3$, 
$
(\Gamma H^*BA)^{W(A)}=\mathbb{F}_5[x_{62}, x_{200}, x_{240}];
$
\item For $p=3$, $G=F_4$, $A=E_{F_4}^3$,
$
(\Gamma H^*BA)^{W(A)}=\mathbb{F}_3[x_{26}, x_{36}, x_{48}];
$
\item For $p=3$, $G=E_6$, $A=E_{3E_6}^{4}$, 
$
(\Gamma H^*BA)^{W(A)}=\mathbb{F}_3[x_{26}, x_{36}, x_{48}, x_{54}];
$
\item For $p=3$, $G=E_7$, $A=E_{2E_7}^{4}$, 
$
(\Gamma H^*BA)^{W(A)}=\mathbb{F}_3[x_{26}, x_{36}, x_{48}, x_{108}];
$
\item For $p=3$, $G=E_8$, $A=E_{E_8}^{5a}$, 
$
(\Gamma H^*BA)^{W(A)}=\mathbb{F}_3[x_{4}, x_{26}, x_{36}, x_{48}, x_{324}],
$
\end{enumerate}
where the subscript of $x$ indicates its cohomological degree.
\end{Proposition}

\begin{proof} We use \fullref{main} for (3), (4) and (5). In these cases,  we described $G_1$, $G_2$, $V_1^*$, $V_2^*$, $\mathbb{F}_{p}[V_1]^{G_1}$, $\mathbb{F}_p[V_2]^{G_2}$ in \fullref{main}.

(1)\qua The case $p=5$, $G=E_8$. $A=\smash{E_{E_8}^3}$. The Weyl group $W(A)$ is the special linear group $SL_{3}(\mathbb{F}_5)$.
The ring of invariants of the special linear group is well-known as Dickson invariants. Then we have
\[
\mathbb{F}_5[t_1, t_2, t_3]^{W(A)} = \mathbb{F}_5[x_{62}, x_{200}, x_{240}],
\]
where  $x_{62}^4=c_{3,0}$, $x_{200}=c_{3,2}$, $x_{240}=c_{3,1}$ and $c_{3, k}$'s are Dickson invariants described in \fullref{odickson}.

(2)\qua The case $p=3$, $G=F_4$, $A=E_{F_4}^3$. The Weyl group $W(A)$ is the special linear group $SL_{3}(\mathbb{F}_3)$.
The ring of invariants are known as Dickson invariants as before:
\[
\mathbb{F}_3[t_1, t_2, t_3]^{W(A)} = \mathbb{F}_3[x_{26}, x_{36}, x_{48}],
\]
where  $x_{26}^2=c_{3,0}$, $x_{36}=c_{3,2}$, $x_{48}=c_{3,1}$ as above.

(3)\qua The case $p=3$, $G=E_6$, $A=\smash{E_{3E_6}^4}$. The Weyl group $W(A)$ is the subgroup of $GL_{4}(\mathbb{F}_3)$
consisting of matrices of the form
\[
\MyMatrix{g_1}{m_0}{g_2},
\]
where $g_1\in G_1=\{1\}$, $g_2\in G_2=SL_{3}(\mathbb{F}_3)$, $m_0\in M_{1,3}(\mathbb{F}_3)$.
Consider $V_1^{*}=\mathbb{F}_3\{t_1\}$ and $V_2^*=\mathbb{F}_3\{t_2, t_3, t_4\}$. Then we have
$\smash{\mathbb{F}_3[V_1]^{G_1}=\mathbb{F}_3[t_1]}$ and $\smash{\mathbb{F}_3[V_2]^{G_2}=\mathbb{F}_3[x_{26}, x_{36}, x_{48}]}$,
where $x_{26}^2=c_{3,0}$, $x_{36}=c_{3,2}$, $x_{48}=c_{3,1}$ and $c_{3,k}$'s are Dickson invariants in $\mathbb{F}_3[V_2]$.
By \fullref{main}, we have
\[
\mathbb{F}_3[t_1, t_2, t_3, t_4]^{W(A)}=\mathbb{F}_3[x_{26}, x_{36}, x_{48}, x_{54}],
\]
where  $x_{54}=\prod_{t\in V_2^*}(t_1+t)$.

(4)\qua The case $p=3$, $G=E_7$, $A=E_{2E_7}^4$. 
The Weyl group $W(A)$ is a subgroup of $GL_{4}(\mathbb{F}_3)$ consisting of matrices of the form
\[
\MyMatrix{g_1}{m_0}{g_2},
\]
where 
$g_1 \in G_1=GL_{1}(\mathbb{F}_3)$, 
$g_2 \in G_2=SL_{3}(\mathbb{F}_3)$ and $m_0\in M_{1,3}(\mathbb{F}_3)$.
Consider $V_1^*=\mathbb{F}_3\{t_1\}$ and  $V_2^*=\mathbb{F}\{t_2, t_3, t_4\}$. Then we have
$\mathbb{F}_3[V_1]^{G_1}=\mathbb{F}_3[t_1^2]$ and $\mathbb{F}_3[V_2]^{G_2}=\mathbb{F}_3[x_{26}, x_{36}, x_{48}]$ as in (3).
By \fullref{main}, we have
\[
\mathbb{F}_3[t_1, t_2, t_3, t_4]^{W(A)}=\mathbb{F}_3[x_{26}, x_{36}, x_{48}, x_{108}],
\]
where $x_{108}=\prod_{t\in V_2^*} (t_1+t)^2$.

(5) The case $p=3$, $G=E_8$, $A=E_{E_8}^{5a}$.
The Weyl group $W(A)$ is a subgroup of $GL_{5}(\mathbb{F}_3)$ consisting of matrices of the form
\[
\MyMatrix{g_1}{m_0}{g_2}=\left( \begin{array}{c|c|c}
g_1 & m_0' & m_o'' \\
\hline
0 & g_2' & 0 \\
\hline
0 & 0 & \epsilon
\end{array}\right),
\]
where $g_1\in G_1=GL_{1}(\mathbb{F}_3)$, 
$g_2 =(g_2', \epsilon) \in G_2=SL_{3}(\mathbb{F}_3)\times GL_{1}(\mathbb{F}_3)\subset GL_{4}(\mathbb{F}_3)$ 
and $m_0=(m_0', m_0'') \in M_{1, 4}(\mathbb{F}_3)=M_{1,3}(\mathbb{F}_3) \times M_{1,1}(\mathbb{F}_3)$.
Consider $V_1^*=\mathbb{F}_3\{ t_1\}$ and $V_2^*=\mathbb{F}_3\{ t_2, t_3, t_4, t_5\}$. Then we have
$
\mathbb{F}_3[V_1]^{G_1}=\mathbb{F}_3[t_1^2]
$
and
$$
\mathbb{F}_3[V_2]^{G_2}=\mathbb{F}_3[x_{26}, x_{36}, x_{48}, t_5^2],
$$
where $x_{26}$, $x_{36}$, $x_{48}$ are Dickson invariants in $\mathbb{F}_3[t_2, t_3, t_4]$ as in (3). By \fullref{main}, we have
\[
\mathbb{F}_3[t_1, t_2, t_3, t_4, t_5]^{W(A)}=\mathbb{F}_3[x_{4}, x_{26}, x_{36}, x_{48}, x_{324}],
\]
where $x_{4}=t_5^2$ and 
$ x_{324}=\prod_{t\in V_2^*} (t_1+t)^2.
$
\end{proof}

\section[Proof of \ref{noncollapsing}]{Proof of \fullref{noncollapsing}}
\label{Sec4}
Let $p$ be a prime, including $p=2$. 
As in the previous section, 
let $G$ be a compact Lie group and $A$ an elementary abelian $p$--subgroup of $G$. 
We denote by $i_{A,G}\co A\to G$ the inclusion.
Then the induced homomorphism \[
\Gamma Bi_{A,G}^{*}\co \Gamma H^{*}\!BG
 \to \Gamma H^{*}\!BA
 \] factors through the ring of invariants of the Weyl group $W(A)$.

\begin{Proposition}  \label{integral} The inclusion of  the image 
 \[
 \mathrm{Im} \Gamma Bi_{A,G}^{*}\to (\Gamma H^{*}\!BA)^{W(A)}
 \] is an integral extension.
 \end{Proposition}
\begin{proof}
By the Peter--Weyl theorem,  for a sufficiently large $n$, there exists an embedding of a compact Lie group $G$
into a unitary group $U(n)$, say
\[
i_{G, U(n)}\co G \to U(n).
\]
Through the induced homomorphism $Bi_{G, U(n)}^{*}\co H^{*}\!BU(n)\to H^*BG$, the mod $p$ cohomology $H^*BG$ is an $H^*BU(n)$--module.
Recall here that
\[
H^{*}\!BU(n)=\mathbb{F}_p[c_1, \ldots, c_n],
\]
where each $c_i$ is a Chern class and $\deg c_i=2i$. So, $H^*BU(n)$ is a Noetherian ring. 
It is well-known that $H^*BG$ is a finitely generated $H^*BU(n)$--module, so that $H^*BG$ is a Noetherian $H^{*}\!BU(n)$--module.
We defined  $\Gamma H^{*}\!BG$ as a quotient module of $H^*BG$.
Therefore, $\Gamma H^*BG$ is also a Noetherian $H^{*}\!BU(n)$--module.
Considering the case $G=A$, we may conclude that 
$\Gamma H^*BA$ is also a Noetherian $H^{*}\!BU(n)$--module.
Since the ring of invariants $(\Gamma H^{*}\!BA)^{W(A)}$ is an $H^*BU(n)$--submodule of a Noetherian $H^*BU(n)$--module $\Gamma H^{*}\!BA$, it is also a Noetherian $H^*BU(n)$--module. 
Hence, the ring of invariants $\smash{(\Gamma H^*BA)^{W(A)}}$ is a finitely generated $H^{*}\!BU(n)$--module.
Thus, the inclusion 
\[
\mathrm{Im}\,\Gamma Bi_{A,G}^* \to (\Gamma H^{*}\!BA)^{W(A)}
\] is an integral extension.
\end{proof}

 In the case $p=3$, $G=E_8$, $\smash{A=E_{E_8}^{5a}}$, the ring of invariants of the Weyl group $W(A)$ is computed in the previous section and 
 it is $$\mathbb{F}_3[x_{4}, x_{26}, x_{36}, x_{48}, x_{324}]$$ as a graded $\mathbb{F}_3$--algebra.
 
Now, we recall the computation of the cotorsion product $$\mathrm{Cotor}_{H^*E_8}(\mathbb{F}_3, \mathbb{F}_3)$$ due to Mimura and Sambe in \cite{MimuraSambe}.
From this, we need just an upper bound for the
 degree of algebra generators of
the cotorsion product. Namely, the following result which is immediate from the computation of Mimura and Sambe suffices to prove the noncollapsing of 
the Rothenberg--Steenrod spectral sequence.

\begin{Proposition}\label{bound}
As a graded $\mathbb{F}_3$--algebra, the cotorsion product $$\mathrm{Cotor}_{H^*E_8}(\mathbb{F}_3, \mathbb{F}_3)$$ 
is generated by elements of degree less than
or equal to $168$.
\end{Proposition}

As a consequence of \fullref{bound}, if the Rothenberg--Steenrod spectral sequence collapsed at the $E_2$--level, 
then  $H^*BE_8$ and $\Gamma H^*BE_8$ 
would be generated by elements of degree less than or equal to $168$ as  graded $\mathbb{F}_3$--algebras. 
The image of the induced homomorphism $\Gamma Bi_{A,E_8}^{*}$ would also be generated by elements of degree less than or equal to $168$.
 Therefore, $\mathrm{Im}\,\Gamma Bi_{A, E_8}^{*}$ would be a subalgebra of $$\mathbb{F}_3[x_{4}, x_{26}, x_{36}, x_{48}].$$
 It is clear that $x_{324}$ is not integral over $\mathbb{F}_3[x_{4}, x_{26}, x_{36}, x_{48}]$,
 and so the inclusion 
 \[
 \mathrm{Im}\,\Gamma Bi_{A, E_8}^{*}\to (\Gamma H^{*}\!BA)^{W(A)}=\mathbb{F}_3[x_{4}, x_{26}, x_{36}, x_{48}, x_{324}]
 \]
would  not be an integral extension. This contradicts \fullref{integral}. 
 Hence, the Rothenberg--Steenrod spectral sequence does not collapse at the $E_2$--level.

\bibliographystyle{gtart}
\bibliography{link}

\begin{thebibliography}{}
\providecommand\bibmarginpar{\leavevmode\marginpar}
\def\urlstyle#1{{\tt #1}}

\bibitem{Andersen}
\textbf{K Andersen}, \textbf{J Grodal}, \textbf{J M\o{}ller}, \textbf{A
  Viruel}, \emph{The classification of $p$--compact groups for $p$ odd}
  \xox{arXiv}{math.AT/0302346}

\bibitem{MimuraSambe}
\textbf{M Mimura}, \textbf{Y Sambe}, \emph{On the cohomology mod $p$ of the
  classifying spaces of the exceptional Lie groups II, III}, J. Math. Kyoto
  Univ. 20 (1980) 327--379 \xox{MR}{582171}

\bibitem{Smith}
\textbf{L Smith}, \emph{Polynomial invariants of finite groups}, Research Notes
  in Mathematics 6, A K Peters Ltd., Wellesley, MA (1995) \xox{MR}{1328644}

\bibitem{Wilkerson}
\textbf{C Wilkerson}, \emph{A primer on the {D}ickson invariants}, from:
  ``Proceedings of the Northwestern Homotopy Theory Conference (Evanston, Ill.,
  1982)'', Contemp. Math. 19, Amer. Math. Soc., Providence, RI (1983)  421--434
  \xox{MR}{711066}

\end{thebibliography}

\end{document}